\documentclass[a4paper,11pt]{article}
\usepackage[latin1]{inputenc}
\usepackage[francais]{babel}
\usepackage{amsmath}
\usepackage{amssymb}
\usepackage{amsfonts}
\usepackage{graphics}
\usepackage{color}
\newcounter{theorem}

\newtheorem{theo}[theorem]{Théorème}%
\newtheorem{prop}[theorem]{Proposition}%
\newtheorem{lemm}[theorem]{Lemme}%
\newtheorem{coro}[theorem]{Corollaire}%
\newtheorem{rema}[theorem]{Remarque}%
\newenvironment{demo}{\trivlist\item[\hskip \labelsep{\rm
D\'emonstration.}\enskip]}%

\def\card{{\rm Card}}

\def\sgn{{\rm sgn}}

\def\ent{{\rm Ent}}

\def\eef{{\bf E}}
\def\iif{{\bf I}}
\def\nnf{{\bf N}}

\def\rrf{{\bf R}}

\def\fc{{\cal F}}

\def\lc{{\cal L}}
\def\mc{{\cal M}}
\def\nc{{\cal N}}

\def\rc{{\cal R}}

\usepackage{graphicx}
\begin{document}
\title{UN PROCESSUS PONCTUEL ASSOCIÉ AUX MAXIMA LOCAUX DU MOUVEMENT BROWNIEN}
\date{}
\author{Christophe LEURIDAN}
\maketitle

\begin{abstract}
Soit $B = (B_t)_{t \in \rrf}$ un mouvement brownien symétrique,
c'est-à-dire un processus tel que $(B_t)_{t \in \rrf_+}$ et
$(B_{-t})_{t \in \rrf_+}$ sont deux mouvements browniens indépendants
issus de $0$. Pour $a \ge b>0$ fixés, nous décrivons la loi de l'ensemble 
$$\mc_{a,b} = \{t \in \rrf : B_t = \max_{s \in [t-a,t+b]} B_s\}.$$ 
Nous relions cet ensemble au fermé régénératif 
$$\rc_a = \{t \in \rrf_+ : B_t = \max_{s \in [(t-a)_+,t]} B_s\},$$
et nous donnons la mesure de Lévy d'un subordinateur dont l'image
fermée est $\rc_a$. 
\end{abstract}

\begin{flushleft}
{\it Classification math\'ematique}~: 60J65, 60G55.\\
{\it Mots-cl\'es}~: mouvement brownien, maximum local, processus
ponctuel, renouvellement, fermé régénératif, subordinateur. 
\end{flushleft}
 
\section{Introduction}

On considère un mouvement brownien symétrique $B = (B_t)_{t \in \rrf}$ :
autrement dit, $(B_t)_{t \in \rrf_+}$ et $(B_{-t})_{t \in \rrf_+}$ sont
deux mouvements browniens indépendants issus de $0$. Le but de cet article est
d'étudier l'ensemble $\mc$ des instants où $(B_t)_{t \in \rrf}$
atteint un maximum local.
On vérifie facilement que presque sûrement, $\mc$ est
dénombrable et dense dans $\rrf$ : il suffit par exemple de montrer
qu'il y a un unique instant réalisant le maximum sur tout intervalle
d'extrémités rationnelles. 

Le fait qu'un instant donné soit un maximum local ne dépend que des
accroissements du mouvement brownien au voisinage de cet
instant. L'indépendance des accroissements du mouvement brownien
suggère donc que l'ensemble $\mc$ est sans mémoire, mais encore
faut-il donner un sens à cette affirmation. Une façon serait de
construire un processus de Poisson ponctuel qui prenne des valeurs
précisément aux instants où $B$ atteint un maximum local. 

Dans \cite{Tsirelson}, Tsirelson montre que cela est possible,
puisqu'on peut construire une suite $(T_n)_{n \ge 1}$ de variables
aléatoires indépendantes, de loi uniforme sur
$[0,1]$ telle que l'ensemble des valeurs $\{T_n ; n \in \nnf\}$ soit
exactement l'ensemble des instants de $[0,1]$ où le mouvement brownien
réalise un minimum local. En fait, cette propriété est partagée par de
nombreux ensembles aléatoires dénombrables denses. Par ailleurs, la
loi uniforme sur $[0,1]$ peut être remplacée par n'importe quelle loi
à densité partout strictement positive sur $[0,1]$. 

Mais la construction de Tsirelson est faite dans un espace
probabilisé plus gros que celui sur lequel est défini le mouvement
brownien : le produit de l'espace de Wiener par $[0,1]^\infty$. 
On peut donc se demander si ce type de construction est possible sans
apport d'aléa extérieur. 
  
Une idée consiste à associer à tout instant $t$ le plus grand
intervalle $I_t$ contenant $t$ sur lequel $B$ reste inférieur ou égal
à $B_t$. Autrement
dit, $I_t = [t-U_t,t+V_t]$, où $U_t = \inf\{s > 0 : B_{t-s} > B_t\}$
et $V_t = \inf\{s > 0 : B_{t+s} > B_t\}$. L'instant $t$ est un maximum
local si et seulement si $\min(U_t,V_t)>0$. Dans ce cas, nous dirons
que $U_t$ et $V_t$ sont les portées à gauche et à droite du maximum
local à l'instant $t$.   

Mais on vérifie facilement que les intervalles $I_s$ et $I_t$ associés à deux
instants $s<t$ sont disjoints si le mouvement brownien dépasse
$\max(B_s,B_t)$ pendant l'intervalle $[s,t]$ et emboîtés dans le cas
contraire. Ce fait est une obstruction au fait que $(\min(U_t,V_t))_{t \in
  \rrf}$ soit un processus de Poisson ponctuel puisqu'il possède une mémoire. 

Dans la suite, nous allons décrire la loi de l'ensemble 
$\mc_{a,b} = \{t \in \rrf : U_t \ge a ; V_t \ge b\}$ pour $a>0$ et $b>0$
fixés. Cet ensemble est formé d'instants isolés, séparés d'au moins
$\min(a,b)$. Nous montrons que les durées entre les instants
successifs de $\mc_{a,b}$ forment des variables aléatoires
indépendantes et de même loi, et nous décrirons cette loi.  

Une méthode pour montrer l'indépendance et l'équidistribution des
durées consiste à relier $\mc_{a,b}$ au fermé régénératif $\rc_a = \{t
\in \rrf_+ : B_t = \max_{s \in [(t-a)_+,t]} B_s\}$. Inversement, nous
utilisons des renseignements obtenus directement sur  $\mc_{a,b}$ pour
décrire la mesure de Lévy d'un subordinateur dont l'image fermée est
le fermé régénératif $\rc_a$. 

Signalons que dans \cite{Neveu-Pitman}, J. Neveu et J. Pitman ont
obtenu des résultats similaires sur les extrema de profondeur
supérieure ou égale à $h$ fixé, appelés « $h$-extrema » par les
auteurs. Avec nos notations, si $t$ est un
instant de maximum local, sa profondeur est la plus petite des quantités 
$$B_t - \min_{s \in [t-U_t,t]} B_s\ \text{ et }\ B_t - \min_{s \in
  [t,t+V_t]} B_s.$$
Pour un minimum local, la définition est symétrique.

\section{Calculs direct de certaines lois associées aux maxima
  locaux.} 

Pour tous réels $s<t$, notons $\rho(s,t)$ le plus petit instant de
$[s,t]$ réalisant le maximum de $B$ sur $[s,t]$ et introduisons les
dénivellations à gauche et à droite correspondantes : 
$$G_{s,t} = B_{\rho(s,t)} - B_s\ \text{ et }\ D_{s,t} = B_{\rho(s,t)} - B_t$$
Ces variables aléatoires sont mesurables pour la tribu $\fc_{s,t}$
engendrée par les accroissements du mouvement brownien entre $s$ et
$t$. Leur utilisation est rendue commode du fait que les tribus
$\fc_{s,t}$ associées à des intervalles d'intérieurs deux à deux
disjoints sont indépendantes. 

L'indépendance permet de montrer facilement le résultat classique
ci-dessous, qui constitue le lemme 13.15 de \cite{Kallenberg}. 

\begin{lemm}\label{Comparaison des maxima locaux} {\bf (Comparaison
    des maxima locaux)}
Preque sûrement, les valeurs des maxima locaux de $B$ sont toutes
différentes. 
\end{lemm}

Comme presque sûrement, le maximum de $B$ sur un segment fixé non
réduit à un point n'est pas atteint aux extrémités, il est donc
réalisé en un unique point intérieur.
 
Le résultat ci-dessous constitue la brique élémentaire des calculs
ultérieurs. 

\begin{prop}
Soient $s<t$ deux réels. Le triplet $(\rho(s,t),G_{s,t},D_{s,t})$ a
même loi que  
$$(s \cos^2 \Theta + t \sin^2 \Theta,\ \sqrt{2X} \sqrt{t-s} \sin \Theta
,\ \sqrt{2Y} \sqrt{t-s} \cos \Theta)$$ 
où $\Theta$, $X$ et $Y$ sont des variables
aléatoires indépendantes, $\Theta$ de loi uniforme sur $[0,\pi/2]$,
$X$ et $Y$ de loi exponentielle de paramètre 1. 
\end{prop}

\begin{demo}
Pour démontrer ce résultat, on peut se limiter au cas où $s=0$ par
stationnarité des accroissements. Mais en notant $S_t = \max\{B_s\ ; s \in
[0,t]\}$, on a $G_{0,t} = S_t$ et $D_{0,t} = S_t-B_t$. 
La densité du triplet $(\rho(0,t),S_t,B_t)$
est connue (voir \cite{Levy} ou le lemme 4 de \cite{Louchard} ou la
généralisation dans \cite{Csaki - Foldes - Salminen}) : en tout
$(r,a,b)$ tel que $0 < r < t$ et $a > \max(0,b)$, elle vaut 
$$\frac{a(a-b)}{\pi r^{3/2} (t-r)^{3/2}} \exp\Big(\frac{-a^2}{2r}\Big)
\exp\Big(\frac{-(a-b)^2}{2(t-r)}\Big).$$
Cela montre que $\rho(0,t)$ suit la loi arcsinus sur l'intervalle
$[0,t]$ et que conditionnellement à $\rho_t=r$, les variables
aléatoires $\frac{S_t^2}{2r}$ et $\frac{(S_t-B_t)^2}{2(t-r)}$ sont
indépendantes de loi exponentielle de paramètre 1. \hfill $\square$
\end{demo}

\begin{coro}
Soient $s<t$ deux réels. Quels que soient $\alpha,\beta \ge 0$,
$$\eef \Big[ \exp \Big(-\frac{\alpha G_{s,t}^2 + \beta D_{s,t}^2}{2}
\Big) \Big] = \frac{1}{\alpha+\beta+\alpha\beta}
\Big(\frac{\alpha}{\sqrt{1+\alpha}} + \frac{\beta}{\sqrt{1+\beta}}
\Big).$$ 
\end{coro}

\begin{demo}
Utilisons l'identité en loi de la proposition. En conditionnant par
rapport à $\Theta$, et en effectuant le changement de variable $t =
\tan \theta$, on obtient 
\begin{eqnarray*}
\eef \Big[\exp \Big(-\frac{\alpha G_{0,1}^2 + \beta D_{0,1}^2}{2} \Big)\Big] 
&=& \eef \big[\exp(-\alpha X \sin^2 \Theta)\ \exp(-\beta Y \cos^2
\Theta) \big]\\
&=& \eef \Big[\frac{1}{1 + \alpha \sin^2 \Theta} \times \frac{1}{1 + \beta 
  \cos^2 \Theta} \Big]\\
&=& \frac{2}{\pi} \int_0^{\pi/2} \frac{d\theta}{(1 + \alpha \sin^2
  \theta)(1 + \beta  \cos^2 \theta)}\\
&=& \frac{2}{\pi} \int_0^\infty 
\frac{dt}{(1+t^2) 
(1 + \frac{\alpha t^2}{1+t^2}) (1 + \frac{\beta}{1+t^2})}\\
&=& \frac{2}{\pi} \int_0^\infty 
\frac{(1+t^2) dt}{\big( 1 + (1+\alpha) t^2 \big) \big((1 + \beta) +
  t^2 \big)}
\end{eqnarray*}
On remarque ensuite que 
$$\frac{\alpha}{1 + (1+\alpha) t^2} + \frac{\beta}{(1 + \beta) + t^2}
= \frac{(\alpha+\beta+\alpha\beta)(1+t^2)}{\big( 1 + (1+\alpha) t^2
  \big) \big(1 + \beta + t^2 \big)},$$
pour en déduire la formule ci-dessus. \hfill $\square$
\end{demo}

\section{Application à l'étude de $\mc_{a,b}$ pour $a \ge b >0$ :
  fonctions de corrélation}

Soient $a \ge b > 0$. Nous allons décrire la mesure de comptage
$N_{a,b}$ associée à $\mc_{a,b}$ : pour tout borélien $A$ de $\rrf$,
$N_{a,b}(A) = \card(\mc_{a,b} \cap A)$. Commençons par quelques
remarques simples.  

\begin{rema} {\bf (Propriétés simples de $N_{a,b}$)}
\begin{enumerate}
\item Deux points de $\mc_{a,b}$ ne peuvent pas être à une distance
  inférieure à $b$. Par conséquent, si le diamètre de $A$ est
  inférieur à $b$, $N_{a,b}(A)$ vaut $0$ ou $1$, suivant que le
  maximum de $B$ sur $A$ ait ou non des portées supérieures ou égales
  à $a$ et $b$. 
\item Comme la loi de $\mc_{a,b}$ est invariante par translation, la
  mesure $A \mapsto \eef[N_{a,b}(A)]$ est un multiple de la mesure de
  Lebesgue. 
\item La variable aléatoire $N_{a,b}(A)$ ne dépend que des
  accroissements du mouvement brownien sur l'ensemble $A + [-a,b]$. A
  cause de l'indépendance des accroissements du mouvement brownien,
  les variables aléatoires comptant les points de $\mc_{a,b}$ dans des
  parties distantes d'au moins $a+b$ sont indépendantes. 
\end{enumerate}
\end{rema}

Nous allons démontrer le résultat suivant 

\begin{theo}~\label{fonction de correlation}
Le processus ponctuel $\mc_{a,b}$ possède des fonctions de corrélation
: pour tout entier $n \ge 1$, 
$$\eef[N_{a,b}(dt_1) \cdots N_{a,b}(dt_n)] =
f_n(t_1,\ldots,t_n)\ dt_1 \cdots dt_n,$$
où $f_n : \rrf^n \to \rrf$ est la fonction symétrique de $n$ variables
réelles définie par 
$$f_n(t_1,\ldots,t_n) = \frac{1}{\pi\sqrt{ab}}\ \prod_{k=1}^{n-1}
h_{a,b}(t_{k+1}-t_k).$$
pour $(t_1,\ldots,t_n)$ tel que $t_1 \le \ldots \le t_n$, avec
$$h_{a,b}(r) = \left|\begin{array}{ccc} 
0 &{\rm si}& r \le b\\
\frac{1}{\pi r} \sqrt{\frac{r-b}{b}} &{\rm si}& b \le r \le a\\
\frac{1}{\pi r} \Big(\sqrt{\frac{r-b}{b}}+\sqrt{\frac{r-a}{a}} \Big) &{\rm si}&
a \le r \le a+b\\
\frac{1}{\pi\sqrt{ab}} &{\rm si} & r \ge a+b
\end{array} \right.$$
\end{theo}

\begin{demo}
Soit $t \in \rrf$. Pour tout $\epsilon \in ]0,b[$, l'intervalle
$[t,t+\epsilon]$ contient au plus un point de $\mc_{a,b}$ et  
$$\rho(t-a,t+b+\epsilon) \in [t,t+\epsilon] \Longrightarrow
N_{a,b}([t,t+\epsilon])=1 \Longrightarrow \rho(t-a+\epsilon,t+b) \in
  [t,t+\epsilon].$$
Ces implications fournissent un encadrement de
$N_{a,b}([t,t+\epsilon])$, 
mais comme cet encadrement est lourd à manipuler,
nous écrirons de façon heuristique 
$$N_{a,b}(dt)=1 \Longleftrightarrow \rho(t-a,t+b) \in dt,$$
d'où 
$$\eef[N_{a,b}(dt)] = \frac{dt}{\pi\sqrt{ab}},$$
ce qui donne la fonction de corrélation pour $n=1$. 

Intéressons-nous à présent au cas où $n \ge 2$. Soient $t_1<\ldots<t_n$. 

Si $t_{k+1}-t_k < b$ pour un certain $k \in [1 \ldots n-1]$, alors
$N_{a,b}(dt_k)N_{a,b}(dt_{k+1}) = 0$ et 
$h_{a,b}(t_{k+1}-t_k)=0$, et le résultat est évident. 

Si $t_{k+1}-t_k > a+b$ pour un certain $k \in [1 \ldots n-1]$, 
alors 
$N_{a,b}(dt_1) \cdots N_{a,b}(dt_k)$ est indépendante de
$N_{a,b}(dt_{k+1}) \cdots N_{a,b}(dt_n)$ d'où  
$$\eef[N_{a,b}(dt_1) \cdots N_{a,b}(dt_n)] =\eef[N_{a,b}(dt_1) \cdots
N_{a,b}(dt_k)] \eef[N_{a,b}(dt_{k+1}) \cdots N_{a,b}(dt_n)].$$
Mais on a aussi $$f_n(t_1,\ldots,t_n) =
f_k(t_1,\ldots,t_k) f_{n-k}(t_{k+1},\ldots,t_n)$$    
grâce au fait que $h_{a,b}(t_{k+1}-t_k)$ est égal à l'intensité
$\frac{1}{\pi\sqrt{ab}}$ du processus ponctuel. 

Ce raisonnement montre qu'on peut donc se limiter au cas où
$b \le t_{k+1}-t_k \le a+b$ pour tout $k \in [1 \ldots n-1]$. Pour
alléger l'écriture, nous nous limitons à deux instants $s<t$ tels
que la différence $r=t-s$, le cas général ne diffèrant que par la
lourdeur des expressions. 

Si $b \le t-s \le a$, alors $s-a \le t-a \le s \le s+b \le t \le
t+b$. Donc 
\begin{eqnarray*} 
\eef[N_{a,b}(ds)N_{a,b}(dt)] 
&=& P[\rho(s-a,s+b) \in ds\ ;\ \rho(t-a,t+b) \in dt]\\
&=& P[\rho(s-a,s+b) \in ds\ ;\ D_{s-a,s+b} < G_{s+b,t+b}\ ;\\
& &\phantom{----} \rho(s+b,t+b) \in dt]
\end{eqnarray*}
Conditionnellement à $(\rho(s-a,s+b),\rho(s+b,t+b))=(s,t)$ les
variables aléatoires $Y = D_{s-a,s+b}^2/2b$ et $X =
G_{s+b,t+b}^2/2(r-b)$ sont indépendantes et de loi exponentielle de paramètre 1 et la probabilité pour que $(r-b)X > bY$ vaut donc $(r-b)/r$. Donc 
$$\eef[N_{a,b}(ds)N_{a,b}(dt)] = \frac{ds}{\pi\sqrt{ab}} \ 
\frac{dt}{\pi\sqrt{(r-b)b}} \ \frac{r-b}{r} = \frac{ds}{\pi\sqrt{ab}}
\ \frac{1}{\pi r} \sqrt{\frac{r-b}{b}} dt.$$

Si $a \le t-s \le a+b$, alors $s-a \le s \le t-a \le s+b \le t \le t+b$. 
En conditionnant par rapport à $\fc_{t-a,s+b}$, on trouve 
\begin{eqnarray*} \eef[N_{a,b}(ds)N_{a,b}(dt)] 
&=& P[\rho(s-a,s+b) \in ds\ ;\ \rho(t-a,t+b) \in dt]\\
&=& P[\rho(s-a,t-a) \in ds\ ;\ D_{s-a,t-a} > G_{t-a,s+b}\ ;\\
& & \phantom{----} D_{t-a,s+b} < G_{s+b,t+b}\ ;\ \rho(s+b,t+b) \in dt]\\
&=& \frac{ds}{\pi\sqrt{a(r-a)}} \times \frac{dt}{\pi\sqrt{b(r-b)}}\\ 
& & \phantom{----} \eef \Big[ \exp \Big( -\frac{G_{t-a,s+b}^2}{2(r-a)} -
\frac{D_{t-a,s+b}^2}{2(r-b)}  \Big) \Big].
\end{eqnarray*}
Mais par changement d'échelle, 
\begin{eqnarray*}
\eef \Big[ \exp \Big( -\frac{G_{t-a,s+b}^2}{2(r-a)} -
\frac{D_{t-a,s+b}^2}{2(r-b)}  \Big) \Big]
&=& \eef \Big[ \exp \Big(-\frac{\alpha G_{0,1}^2 + \beta D_{0,1}^2}{2}
\Big) \Big] 
\end{eqnarray*}
avec 
$$\alpha = \frac{a+b-r}{r-a}, \phantom{--} \beta = \frac{a+b-r}{r-b}.$$
Comme
$$1+\alpha = \frac{b}{r-a}, \phantom{--}\ 1+\beta = \frac{a}{r-b},$$
$$\alpha + \beta + \alpha\beta = (1+\alpha)(1+\beta)-1 =
\frac{ab-(r-a)(r-b)}{(r-a)(r-b)} =
\frac{r(a+b-r)}{(r-a)(r-b)},$$
on a donc d'après le corollaire
\begin{eqnarray*}
\eef \Big[ \exp \Big(-\frac{\alpha G_{0,1}^2 + \beta D_{0,1}^2}{2}
\Big) \Big] &=& \frac{(r-a)(r-b)}{r(a+b-r)}\ \Big(
\frac{a+b-r}{\sqrt{b(r-a)}} + \frac{a+b-r}{\sqrt{a(r-b)}} \Big)\\
&=& \frac{\sqrt{(r-a)(r-b)}}{r\sqrt{ab}} (\sqrt{a(r-b)} + \sqrt{b(r-a)}),
\end{eqnarray*}
d'où le résultat. \hfill $\square$  
\end{demo}

\begin{rema}{\bf (Etude du maximum de la fonction $h_{a,b}$)}

Si $a \le 4b$, la fonction $h_{a,b}$ est majorée par sa limite en
$+\infty$. Autrement dit, la probabilité pour qu'il y ait un point
de $\mc_{a,b}$ au voisinage d'un instant $t$ sachant que $s \in
\mc_{a,b}$ est maximum pour $|t-s| \ge a+b$ : les points de
$\mc_{a,b}$ « se repoussent ». 

En revanche, si $a>4b$, la fonction $h_{a,b}$ possède un maximum
global strict en $2b$ : la probabilité pour qu'il y ait un point
de $\mc_{a,b}$ au voisinage d'un instant $t$ sachant que $s \in
\mc_{a,b}$ est maximum pour $|t-s| = 2b$.
\end{rema}

\begin{demo}
On vérifie facilement que la fonction $r \mapsto \sqrt{a(r-b)} +
\sqrt{b(r-a)}$ est concave sur $[a,+\infty[$, que sa valeur en $a+b$
est $a+b$ et que sa dérivée en $a+b$ vaut $1$. Par conséquent, pour
tout $r \in [a,+\infty[$, $\sqrt{a(r-b)} + \sqrt{b(r-a)} \le r$, d'où
$h_{a,b}(r) \le 1/(\pi \sqrt{ab})$. 

Par ailleurs, on vérifie facilement que la fonction $r
\mapsto r^{-1}\sqrt{(r-b)}$ est strictement croissante sur $[b,2b]$,
strictement décroissante sur $[2b,+\infty[$, si bien que le maximum de
$\pi \sqrt{b} h_{a,b}$ sur $[b,a]$ est majoré par $1/2\sqrt{b}$, avec
égalité lorsque $a \ge 2b$. Pour que ce maximum excède $1/\sqrt{a}$,
il faut et il suffit que $a > 4b$. \hfill $\square$
\end{demo}

On peut se demander si le processus ponctuel $\mc_{a,b}$ est
déterminantal, comme par exemple l'ensemble
des zéros d'une série entière à coefficients indépendants de loi
$\nc(0,1)$ étudié dans \cite{Peres-Virag}. La réponse est négative. 

\begin{rema}
Le processus ponctuel $\mc_{a,b}$ n'est pas déterminantal :
on ne peut pas trouver d'application symétrique $K$ de $\rrf^2$ dans
$\rrf$ tel que les fonctions de corrélation soient données par   
$f_n(t_1,\ldots,t_n) = \det \big(K(t_i,t_j)\big)_{1 \le i,j \le n}$.
\end{rema}

\begin{demo}
En effet, si une telle application $K$ existait, elle devrait vérifier
pour tout $t \in \rrf$
$$K(t,t) = f_1(t) = c\ \text{ avec }\ c = \frac{1}{\pi\sqrt{ab}}$$
et pour tout $s,t \in \rrf$
$$c^2 - K(s,t)^2 = f_2(s,t)$$
Mais $f_n(t_1,\ldots,t_n)>0$ si et seulement si $t_1,\ldots,t_n$ sont
à distance $>b$ les uns des autres. 
On aurait donc $K(0,b) = \xi c$, $K(b,2b) = \eta c$, $K(0,2b) = \zeta
c$ avec $|\xi|=|\eta|=1$ et $|\zeta|<1$, d'où 
$$f_3(0,b,2b) = c^3 \left|\begin{array}{ccc}
1    &\xi &\zeta\\
\xi  &1   &\eta\\
\zeta&\eta&1
\end{array} \right| = -c^3 (1 - 2\xi\eta\zeta + \zeta^2) < 0,$$
ce qui est absurde. \hfill $\square$
\end{demo}

\section{Indépendance et loi des durées entre les points de
  $\mc_{a,b}$}

La forme des fonctions de corrélation, et plus précisément le fait 
que la fonction de corrélation à $n$ points se factorise en 
$f_n(t_1,\ldots,t_n) = c\ h(t_2-t_1) \cdots h(t_n-t_{n-1})$ pour $t_1
< \ldots < t_n$ permet de montrer que les durées successives entre les
instants de $\mc_{a,b}$ sont indépendantes et de même loi. Ces
instants forment un processus de renouvellement stationnaire, comme le
processus des $h$-extrema introduit par J. Neveu et J. Pitman dans
\cite{Neveu-Pitman}. Mais la loi des durées entre les instants
successifs de $\mc_{a,b}$ est plus compliquée que la loi des durées
entre les $h$-extrema successifs, qui a pour transformée de Laplace
$\theta \mapsto 1/{\rm ch}(h\sqrt{2 \theta})$. 

%Dans la
%partie~\ref{lien entre R(a) et M(a,b)}, nous retrouverons ce
%résultat en reliant l'ensemble $\mc_{a,b}$ à un subordinateur pour
%utiliser le caractère poissonien des sauts. 

\begin{theo}\label{Loi des durees}{\bf (Loi de $\mc_{a,b}$)}\\ 
Notons $\ldots<T_{-1}<T_0<T_1<T_2<\ldots$ les instants successifs
de $\mc_{a,b}$ avec $T_0 \le 0 < T_1$. Alors 
\begin{itemize}
\item Les durées $\ldots,T_{-1}-T_{-2}, T_0-T_{-1},
  T_2-T_1, T_3-T_2, \ldots$ sont indépendantes, de même
  loi et forment une suite indépendante du couple $(T_0,T_1)$. 
\item La durée $T_2-T_1$ possède une densité donnée par 
$$g_{a,b}(r) = \sum_{n=1}^\infty (-1)^{n-1} h_{a,b}^{*n}(r)$$ où
$h_{a,b}$ est la fonction donnée dans le théorème~\ref{fonction de
  correlation}.
\item La variable $(T_0,T_1)$ possède une
densité donnée par 
$$f_{a,b}(t_0,t_1) = \frac{\iif_{t_0<0<t_1}}{\pi\sqrt{ab}}
g_{a,b}(t_1-t_0).$$ 
\end{itemize}
\end{theo}

\begin{demo}
Dans toute la démonstration, nous omettrons les indices $a,b$ et nous
poserons $c = 1/(\pi\sqrt{ab})$. La démonstration comprend deux
étapes. 

1. Commençons par obtenir la loi d'un $n$-uplet $(T_1,\ldots,T_n)$. 
Le passage des fonctions de corrélation à la loi de
$(T_1,\ldots,T_n)$ se fait par des arguments classiques figurant par
exemple dans~\cite{Daley - Vere-Jones} au
chapitre 5, section 4 « moment measures and product densities »,
notament l'exemple 5.4(b). Voyons ci-dessous une démonstration
directe, inspirée d'idées qui m'ont été suggérées par J. Brossard. 

Pour $n \in \nnf^*$, notons $D_n = \{(t_1,\ldots,t_n) \in \rrf^n :
t_1<\ldots<t_n\}$. A l'aide du théorème des classes monotones, on
vérifie que pour tout borélien $C \subset D_n$,  
$$\eef[\card((\mc_{a,b})^n \cap C)] = \int_C f_n(t_1,\ldots,t_n)\ dt_1
\cdots dt_n.$$
Calculons la loi de $(T_1,T_2)$ ; un calcul semblable fournit la loi de
$(T_1,\ldots,T_n)$, la seule différence étant la longueur des
formules.  

Soient $I_1 = [a_1,b_1]$ et $I_2 = [a_2,b_2]$ des intervalles de longueur $<
b$ dont les bornes vérifient $0 < a_1 < b_1 < a_2 < b_2$. Notons $J_1
= ]0,a_1[$ et $J_2 = ]b_1,a_2[$.   
Comme $I_1$ et $I_2$ contiennent au plus un point de $\mc_{a,b}$, on a 
$$\iif_{[T_1 \in I_1\ ; T_2 \in I_2]} 
= \iif_{[N(J_1)=0\ ;\ N(I_1)=1\ ;\ N(J_2)=0\ ;\ N(I_2)=1]}.$$
En utilisant le fait que $0^0=1$, $0^m=0$ pour $m \ge 1$, et la
formule du binôme, on obtient : 
\begin{eqnarray*}
\iif_{[T_1 \in I_1\ ; T_2 \in I_2]} 
&=& (1-1)^{N(J_1)} N(I_1) (1-1)^{N(J_2)} N(I_2)\\
&=& \sum_{k,l \in \nnf} (-1)^{k+l} \left({N(J_1) \atop
    k} \right) N(I_1)  \left({N(J_2)\atop l} \right) N(I_2)\\  
&=& \sum_{k,l \in \nnf} (-1)^{k+l} \card((\mc_{a,b})^{k+l+2} \cap C_{k,l}),
\end{eqnarray*}
avec $C_{k,l} = D_{k+l+2} \cap J_1^k \times I_1 \times J_2^l \times
I_2$.
En passant aux espérances, on obtient donc 
$$P[T_1 \in I_1\ ; T_2 \in I_2] = \sum_{k,l \in \nnf} (-1)^{k+l} e_{k,l}$$ 
avec 
\begin{eqnarray*}
e_{k,l} &=&\eef[\card((\mc_{a,b})^{k+l+2} \cap C_{k,l}]\\
&=& \int_{C_{k,l}} f_{k+l+2} (u_1,\ldots,u_k,t_1,v_1,\ldots,v_l,t_2)
\ du_1 \cdots du_k\ dt_1\ dv_1 \cdots dv_l\ dt_2\\
&=& \int_{I_1 \times I_2} \big( \int_{0<u_1<\ldots<u_k<a_1} 
c\ h(u_2-u_1) \cdots h(u_k-u_{k-1}) h(t_1-u_k)\ du_1 \cdots du_k \big)\\ 
& &\!\!\!\!\big( \int_{b_1<v_1<\ldots<v_l<a_2} h(v_1-t_1)
h(v_2-v_1) \cdots h(v_l-v_{l-1}) h(t_2-v_l)\ dv_1 \cdots dv_l \big) dt_1 dt_2. 
\end{eqnarray*}
Comme le support de $h$ est contenu dans $[b,+\infty[$ et comme les
intervalles $I_1$ et $I_2$ sont de longueur $<b$, on ne change pas les
intégrales ci-dessus en remplaçant les domaines d'intégration donnés
par les inégalités
$$0<u_1<\ldots<u_k<a_1\ \text{ et }\ b_1<v_1<\ldots<v_l<a_2$$
par les domaines légèrement plus gros définis par 
$$0<u_1<\ldots<u_k<t_1\ \text{ et }\ t_1<v_1<\ldots<v_l<t_2.$$
Les intégrales obtenues s'interprètent comme des produits de
convolution et on trouve 
\begin{eqnarray*}
e_{k,l} &=& \int_{I_1 \times I_2} (c\ \iif_{\rrf_+} * h^{*k}) (t_1)
\ h^{*(l+1)}(t_2-t_1)\ dt_1\ dt_2. 
\end{eqnarray*}
Ainsi, 
$$P[T_1 \in I_1\ ; T_2 \in I_2] = \int_{I_1 \times I_2} \Big(\sum_{k \in
  \nnf} (-1)^k\ c\ \iif_{\rrf_+} * h^{*k} (t_1) \Big) \Big(\sum_{l \in
  \nnf} (-1)^l\ h^{*(l+1)} (t_2-t_1)\Big)\ dt_1\ dt_2,$$ 
ce qui montre que les variables aléatoires $T_1$ et $T_2-T_1$ sont
indépendantes, que $T_2-T_1$ admet comme densité
$$g = \sum_{l \in \nnf} (-1)^l\ h^{*(l+1)} = \sum_{k \in \nnf^*}
(-1)^{k-1}\ h^{*k}$$ 
et que $T_1$ admet comme densité
$$f = \sum_{k \in \nnf} (-1)^k\ c\ \iif_{\rrf_+} * h^{*k} = c\ (\iif_{\rrf_+}
- \iif_{\rrf_+} * g).$$
Remarquons que les interversions des sommes avec les espérances ou
avec les intégrales ne posent pas de problème : toutes les sommes se
ramènent à des sommes finies du fait que $N(J) \le \ent(L/b)+1$ si $J$
est un intervalle de longueur $L$ et grâce au fait que $h$ est à
support dans $[b,+\infty[$. 

2. Déterminons maintenant la loi du $2n$-uplet
$(T_{-n+1},\ldots,T_0,T_1,\ldots,T_n)$.

Soient $I_{-n+1},\ldots,I_0,I_1,\ldots,I_n$ des intervalles de
$\rrf$, non enchevétrés, rangés dans cet ordre tels que $0$ majore
$I_0$ et minore $I_1$. Notons $m$ la borne inférieure de $I_{-n+1}$ et
supposons que $I_{-n+1}$ est de longueur $<b$. Dans ce cas, $I_{-n+1}$ contient
au plus un point de $\mc_{a,b}$, si bien que l'événement 
$$[\forall k \in [-n+1 \ldots n], T_k \in I_k].$$
signifie que $T_{-n+1},\ldots,T_n$ sont les instants
successifs de $\mc_{a,b}$ après $m$. Par stationnarité de $\mc_{a,b}$,
cet événement a pour probabilité   
$$P[\forall k \in [1 \ldots 2n], T_k \in J_k] = \int_{J_1 \times \cdots
  \times J_{2n}} f_{}(s_1)g_{}(s_2-s_1)\cdots g_{}
(s_{2n}-s_{2n-1}) ds_1 \cdots ds_{2n}$$
avec $J_k = I_{k-n}-m$ pour tout $k \in [1 \ldots 2n]$. Mais $J_1
\subset [0,b]$, d'où $f_{}(s_1) = c$ pour tout
$s_1 \in J_1$. En effectuant le changement de variables $s_k =
t_{k-n}-m$, on obtient donc
$$P[\forall k \in [-n+1 \ldots n], T_k \in I_k] = \int_{I_{-n+1} \times \cdots
  \times I_{n}} c \prod_{|k| \le n-1} g_{}(t_{k+1}-t_k) dt_{-n+1} \cdots dt_{n}.$$
Comme $g_{}$ est nulle sur $\rrf_-$, on en déduit que
$(T_{-n+1},\ldots,T_n)$ admet pour densité 
$$(t_{-n+1},\ldots,t_n) \mapsto \frac{\iif_{t_0<0<t_1}}{\pi\sqrt{ab}}
\prod_{|k| \le n-1} g_{}(t_{k+1}-t_k),$$
d'où le résultat. \hfill $\square$ 
\end{demo}

\paragraph{Remarque.} L'invariance de la loi de $\mc_{a,b}$ par
translation permet de déduire facilement la densité de $(T_0,T_1)$ de
celle de $T_1$. En effet, pour $t_0 < 0 < t_1$, 
$$P[T_0 \le t_0 ; T_1 > t_1] = P[\mc_{a,b} \cap ]t_0,t_1] = \emptyset] =
P[\mc_{a,b} \cap ]0,t_1-t_0] = \emptyset] = P[T_1 > t_1-t_0].$$

Dans la suite, nous allons démontrer autrement l'indépendance des
durées entre les points successifs de $\mc_{a,b}$. Pour cela, nous
allons nous intéresser à l'ensemble $\mc_{a,0} =
\{t \in \rrf : U_t \ge a\}$ des instants qui sont des records depuis
une durée au moins égale à $a$. L'étude de cet ensemble qu'on pourrait
qualifier de « fermé stationnaire régénératif » se ramène à celle du
fermé régénératif $\rc_a = \{t \in \rrf_+ : U_t \ge \min(t,a)\}$.

\section{Le fermé régénératif $\rc_a$}

Dans cette partie, nous notons $A^a_t = \max\{B_s\ ; s \in
[(t-a)_+,t]\}$ pour tout $t \in \rrf_+$, si bien que $\rc_a = \{t \in
\rrf_+ : A^a_t = B_t\}$. Pour ne pas alourdir les notations, nous
omettrons souvent l'indice $a$ lorsqu'il n'y a pas d'ambiguïté. 

\begin{prop} 
Le fermé $\rc_a$ est régénératif dans la filtration naturelle de
$(B_t)_{t \in \rrf_+}$, notée $(\fc^B_t)_{t \in \rrf_+}$. 
\end{prop}

\begin{demo}
Comme $\rc_a$ est l'ensemble des zéros de $(A_t-B_t)_{t \in \rrf_+}$,
on en déduit que $\rc_a$ est prévisible (donc progressif) dans la
filtration naturelle de $B$.

Soient $r>0$ et $D_r = \inf(\rc_a \cap ]r,+\infty[)$. 
Notons $B'_t = B_{D_r+t} - B_{D_r}$, $A'_t = \max\{B'_s\ ; s \in
[(t-a)_+,t]\}$ pour tout $t \in \rrf_+$ et $\rc'_a = \{t \in \rrf_+ :
A'_t = B'_t\}$. Alors pour tout $r \ge 0$, 
$$D_r+t \in  \rc_a 
\Longleftrightarrow B_{D_r+t} = \max\{B_s\ ; s \in [(D_r+t-a)_+,D_r+t]\}$$
Mais le maximum de $B$ sur l'intervalle $[(D_r+t-a)_+,D_r+t]$ est
aussi le maximum de $B$ sur l'intervalle $[D_r+(t-a)_+,D_r+t]$.
En effet, si $t \le a$, alors $(D_r - a)_+ \le (D_r+t-a)_+ \le D_r$
donc $B_{D_r} = \max\{B_s\ ; s \in [(D_r+t-a)_+,D_r]\}$ ; si $t \ge
a$, il n'y a rien à prouver. Donc 
$$D_r+t \in  \rc_a 
\Longleftrightarrow B'_t = \max\{B'_s\ ; s \in [(t-a)_+,+t]\}$$
ce qui montre que $(\rc_a - D_r)_+ = \rc'_a$ est indépendant de
$\fc^B_{D_r}$ et de même loi que $\rc_a$. \hfill $\square$
\end{demo}

\begin{prop} (Propriétés du processus $(A^a_t)_{t \in \rrf_+}$)
\begin{enumerate}
\item Le processus $(A^a_t)_{t \in \rrf_+}$ est continu et à variation
  localement bornée. 
\item La décomposition $(A^a_t)_{t \in \rrf_+}$ en somme d'un processus
croissant et d'un processus décroissant s'écrit 
$$A^a_t = \int_0^t \iif_{[A^a_s=B_s]} dA^a_s + \int_0^t
\iif_{[A^a_s=B_{s-a}]} dA^a_s.$$ 
\item Le processus $L^a = \int_0^\cdot \iif_{[A^a_s=B_s]} dA^a_s$ est
  le temps local symétrique en $0$ de la semi-martingale
  $(A^a_t-B_t)_{t \in \rrf_+}$. 
\end{enumerate}
\end{prop}

\begin{demo} 
La démonstration comporte plusieurs étapes.

La continuité du processus $(A_t)_{t \in \rrf_+}$ est évidente. 

L'ensemble des instants $t$ tels que $A_t > B_t$ est un ouvert de
$\rrf_+$ et le processus $(A_t)_{t \in \rrf_+}$ est décroissant sur
chaque composante connexe de cet ouvert. En effet, si $A_{t_0} >
B_{t_0}$, alors par continuité de $B$, $A_t =
\max\{B_s\ ;\ s \in [(t-a)_+,t_0]\}$ au voisinage de $t_0$. 
De même, l'ensemble des instants $t$ tels que $A_t >
  B_{(t-a)_+}$ est un ouvert de $\rrf_+$ et le processus $(A_t)_{t \in
    \rrf_+}$ est croissant sur chaque composante connexe de cet ouvert. 

Définissons une suite croissante de temps d'arrêt
  $(T_n)_{n \in \nnf}$ par $T_0=0$ et 
$$T_{2n+1} = \inf\{t > T_{2n} : A_t = B_{(t-a)_+}\}.$$
$$T_{2n+2} = \inf\{t \ge T_{2n+1} : A_t = B_t\}.$$ 
D'après les remarques précédentes, le processus $A$ est croissant
sur chaque intervalle $[T_{2n},T_{2n+1}]$, et décroissant
sur chaque intervalle $[T_{2n+1},T_{2n+2}]$. Par ailleurs, presque sûrement
sur l'événement $[T_{2n} < \infty]$, pour tout $t \in ]T_{2n},T_{2n}+a]$ 
$$A_t \ge \max\{B_s\ ;\ s \in [T_{2n},t]\} > B_{T_{2n}} = A_{T_{2n}}
\ge B_{(t-a)_+}$$
ce qui entraîne que $T_{2n+1} \ge T_{2n}+a$. Par conséquent $T_n \to
+\infty$ quand $n \to +\infty$. Ainsi, le processus $A$ est monotone
par morceaux.  

L'intersection des fermés $\{t \in
\rrf_+ : A_t=B_t\}$ et $\{t \in \rrf_+ : A_t=B_{(t-a)_+}\}$ est au
plus dénombrable : elle est contenue dans l'ensemble des instants
$T_n$ et même dans l'ensemble des instants de la forme $T_{2n}$,
puisque l'égalité $A_{T_{2n+1}} = B_{T_{2n+1}}$ entraîne $T_{2n+2} =
T_{2n+1}$. Comme le processus $(A_t)_{t \in \rrf_+}$ est
constant sur tout intervalle où $A_t$ est différent de $B_t$ et
$B_{t-a}$, il se décompose sous la forme
$$A_t = \int_0^t \iif_{[A_s=B_s]} dA_s + \int_0^t
\iif_{[A_s=B_{s-a}]} dA_s.$$ 

Le temps local symétrique en $0$ de la semi-martingale $(A_t-B_t)_{t
  \in \rrf_+}$ est donné par la formule de Tanaka
$$A_t - B_t = A_0 - B_0 + \int_0^t \sgn(A_s-B_s) d(A_s-B_s) + L_t,$$
avec la convention $\sgn(0)=0$. Comme pour tout $s>0$, $P[A_s=B_s]=0$,
l'ensemble des $s \in \rrf_+$ tels que $A_s=B_s$ est presque sûrement
de mesure nulle. En utilisant la décomposition du processus $(A_t)_{t
  \in \rrf_+}$ démontrée plus haut, on obtient
$$A_t - B_t = \int_0^t \iif_{[A_s=B_{(s-a)_+}]} dA_s - B_t + L_t,$$
d'où par différence, 
$$L_t = \int_0^t \iif_{[A_s=B_s]} dA_s,$$ 
ce qui termine la démonstration. \hfill $\square$
\end{demo}

\begin{rema}
Presque sûrement, l'intersection des fermés $\{t \in \rrf_+ : A^a_t=B_t\}$ et
$\{t \in \rrf_+ : A^a_t=B_{(t-a)_+}\}$ est réduite à $\{0\}$. 
\end{rema}

Bien que ce résultat ne soit pas utile pour la suite, nous
donnons une démonstration, qui utilise la fragmentation
en intervalles associée à l'excursion brownienne, (voir~\cite{Bertoin 2}). 

\begin{demo}
Grâce à la propriété de Markov, il suffit de montrer que presque
sûrement, $T_2$ n'appartient pas à $\{t \in \rrf_+ : A_t=B_{(t-a)_+}\}$. 

Notons $\epsilon$ la première excursion de longueur $\ge a$ du mouvement
brownien réfléchi $(S_t-B_t)_{t \in \rrf_+}$, $\Lambda$ sa longueur et
${\bf e}$ l'excursion brownienne de longueur $1$ obtenue par
changement d'échelle à partir de $\epsilon$. Les variables aléatoires
$\Lambda$ et ${\bf e}$ sont indépendantes, et $P[\Lambda > x] =
\sqrt{a/x}$ pour tout $x \ge a$. En particulier $\Lambda >a$ presque
sûrement. 

L'excursion $\epsilon$ débute à l'instant $T_1-a$ et finit après
l'instant $T_2$. Plus précisément, regardons la fragmentation en
intervalles associée à l'excursion $\epsilon$ : pour tout $h>0$, on
note $F_h$ la collection des intervalles d'excursions de $(S_t-B_t)_{t
  \in \rrf_+}$ au-dessus de $h$ contenues dans l'excursion
$\epsilon$. Alors $T_2$ est la fin du premier intervalle de longueur
$> a$ au moment où celui-ci se casse en deux intervalles de longueur
$\le a$. Pour que $A_{T_2}=B_{T_2-a}$, il faudrait que la
fragmentation associée à l'excursion ${\bf e}$ produise un intervalle
de longueur exactement égale à $a/\Lambda$. 

Mais si l'on choisit $U$ uniformément dans $[0,1]$ et indépendamment
de ${\bf e}$, on sait que la longueur de l'intervalle contenant $U$ au
cours de la fragmentation évolue à un changement de temps près
comme $e^{-\xi}$, où $\xi$ est un subordinateur. Le subordinateur
$\xi$ ne dépend que de ${\bf e}$ et de $U$ et est donc indépendant de
$\Lambda$. Comme il est sans dérive, $P[\exists t
\in \rrf_+ : e^{-\xi_t} = a/\Lambda] = 0$, d'où
$P[B_{T_2}=B_{T_2-a}]=0$.  \hfill $\square$
\end{demo}

On sait que tout fermé régénératif parfait est l'image fermée d'un
subordinateur, unique à changement de temps linéaire près : voir le
théorème 2.1 de \cite{Bertoin}, démontré dans l'article
\cite{Maisonneuve} de B. Maisonneuve. Une façon d'obtenir le
subordinateur est de construire un temps local associé au fermé
régénératif et de prendre son inverse. Nous allons voir que le temps
local $L^a$ fait l'affaire, bien que qu'il soit défini comme temps
local de la semimartingale $A^a-B$. 

\begin{prop} 
L'inverse continu à droite du processus $L^a = \int_0^t
\iif_{[A^a_s=B_s]} dA^a_s$, défini par $\sigma^a_l = \inf\{t \ge 0 :
L^a_t>l\}$ est un subordinateur dont l'image fermée est $\rc_a$. 
\end{prop}

\begin{demo}
L'image fermée de $\sigma^a$ est égale au support de la mesure de
Stieltjes associée à $L$, qui est inclus dans $\rc_a = \{t \in
\rrf_+ : A_t = B_t\}$. A l'aide de la propriété forte de
Markov, on montre que pour tout rationnel $r>0$, l'instant 
$D_r = \inf(\rc_a \cap ]r,+\infty[)$ est (presque sûrement) un instant
de croissance des processus $A$ et $L$, donc appartient à l'image de
$\sigma^a$. Comme $\rc_a$ est d'intérieur vide, tout instant de
$\rc_a$ peut être approché par un instant de la forme $D_r$, si bien
que $\rc_a$ est contenu dans l'image fermée de $\sigma^a$. 

Comme $L_t \ge A_t \ge B_t$ pour tout $t \in \rrf_+$ et comme presque
sûrement, les trajectoires browniennes ne sont pas bornées, le
processus croissant $L$ tend vers $\infty$, si bien que les temps
d'arrêt $\sigma^a_l$ sont finis presque sûrement. Fixons $l>0$ et 
notons pour $t \ge 0$ 
$$B'_t = B_{\sigma^a_l+t} - B_{\sigma^a_l},$$ 
$$A'_t = \max\{B'_s\ ; s \in [(t-a)_+,t]\}$$ 
$$L'_t = \int_0^t \iif_{[A'_s=B'_s]} dA'_s.$$
Comme dans la proposition 1, on montre que pour tout $t \ge 0$,
$A_{\sigma^a_l+t} = B_{\sigma^a_l} + A'_t$ grâce à l'égalité $A_{\sigma^a_l}
= B_{\sigma^a_l}$. On en déduit que
$$L_{\sigma^a_l+t} - A_{\sigma^a_l} = \int_0^t
\iif_{[A_{\sigma^a_l+s}=B_{\sigma^a_l+s}]} dA_{\sigma^a_l+s} = \int_0^t
\iif_{[A'_s=B'_s]} dA'_s = L'_t.$$
Par conséquent, le processus $\sigma^a_{l+\cdot} - \sigma^a_l$ est
l'inverse continu à droite de $L'$. Ce processus est donc 
indépendant de $\fc^B_{\sigma^a_l}$ et a même loi que $\sigma^a$.
\hfill $\square$
\end{demo}

Nous décrirons plus loin la mesure de Lévy du subordinateur
$\sigma^a$. Donnons déjà un résultat immédiat. 

\begin{prop}
Soit $\nu_a$ la mesure de Lévy du subordinateur $\sigma^a$. Pour tout $b
\in ]0,a]$, 
$$\nu_a[b,\infty[ = \sqrt{\frac{2}{\pi b}}.$$ 
\end{prop}

\begin{demo}
On vérifie facilement que le processus $A$ coïncide avec le processus $S$
défini par $S_t = \max\{B_s\ ; s \in [0,t]\}$ jusqu'à l'instant 
$T_1 = \inf\{t \ge a : A_t = B_{t-a}\}$, qui correspond au premier palier de
longueur $\ge a$ des deux processus. Leurs inverses continus à
droite coïncident donc jusqu'au premier saut de taille $\ge a$, qui a lieu au
même moment pour les deux. Donc les mesures de Lévy de ces deux
subordinateurs donnent la même mesure à l'intervalle $[b,+\infty[$
pour tout $b \in ]0,a]$, puisque le premier saut de hauteur $\ge b$
d'un subordinateur de mesure de Lévy $\nu$ se produit au bout d'un
temps exponentiel de paramètre $\nu[b,+\infty[$. \hfill $\square$
\end{demo}

\section{Lien entre $\rc_a$ et  $\mc_{a,b}$ pour $a \ge b >0$}~\label{lien
  entre R(a) et M(a,b)} 

Supposons que $a \ge b >0$. Sur l'événement presque sûr où les maxima
locaux sont à des hauteurs toutes différentes (voir proposition
~\ref{Comparaison des maxima locaux}), on vérifie facilement
l'équivalence suivante, valable pour tout $t>a$ :  
$$t \in \mc_{a,b} \Longleftrightarrow \rc_a
\cap [t,t+b[ = \{t\}.$$
Par conséquent la trace de $\mc_{a,b}$ sur $]a,+\infty[$ est
l'ensemble des débuts des intervalles de longueur $\ge b$ dans
l'ouvert $\rc_a^c\ \cap\ ]a,+\infty[$. 

Ces intervalles correspondent aux sauts de hauteur $\ge b$ du
subordinateur $\sigma^a$ après le franchissement de $a$. Plus
précisément, les éléments de $\mc_{a,b} \cap ]a,+\infty[$ sont
les instants $T'_n = \sigma^a_{(R_n)-}$ pour $n \in \nnf^*$ où
$(R_n)_{n \in \nnf}$ est la suite définie par 
$R_0 = \inf\{l \ge 0 : \sigma^a_l>a\}$ et
$R_n = \inf\{l > R_{n-1} : \Delta\sigma^a_l \ge b\}$. 
Comme le processus des sauts $(\Delta\sigma^a_l)_{l \ge 0}$ est un
processus de Poisson ponctuel, on en déduit que les durées
$T'_2-T'_1,T'_3-T'_2,\ldots$ forment une suite de variables aléatoires
indépendantes et de même loi, indépendante de $T'_1$. 

Soit $\nu_a$ la mesure de Lévy du subordinateur $\sigma^a$. On peut
exprimer la transformée de Laplace de $T'_2-T'_1$ à l'aide $\nu_a$ en écrivant
$$T'_2-T'_1 = \Delta\sigma^a_{R_1} + \lim_{\epsilon \to 0}\ \sum_{R_1<l<R_2}
\Delta\sigma^a_l \iif_{[\Delta\sigma^a_l \ge \epsilon]}.$$
En effet :
\begin{itemize}
\item le saut $\Delta\sigma^a_{R_1}$ a pour loi $\nu_a(\cdot|[b,\infty[)$ ;
\item les sauts successifs de hauteur $\ge \epsilon$ effectués par $\sigma^a$ après l'instant $R_1$ ont pour loi $\nu_a(\cdot|[\epsilon,\infty[)$ ;
\item toutes ces variables aléatoires sont indépendantes. 
\end{itemize}
Le nombre de sauts de hauteur $\ge \epsilon$ jusqu'au premier saut de hauteur $\ge b$ suit une loi géométrique de paramètre
  $\nu_a[b,\infty[/\nu_a[\epsilon,\infty[$. 
Donc pour tout $\theta \ge 0$, $\exp[-\theta(T'_2-T'_1)]$ est la
limite quand $\epsilon \to 0$ de
\begin{eqnarray*} 
\int_{\rrf_+^*} e^{-\theta x}\nu_a(dx|[b,\infty[) &\times& \sum_{n=1}^\infty
\frac{\nu_a[b,\infty[}{\nu_a[\epsilon,\infty[}
\Big(\frac{\nu_a[\epsilon,b\ [}{\nu_a[\epsilon,\infty[}\Big)^{n-1} \Big(
\int_{\rrf_+^*} e^{-\theta x}\nu_a(dx|[\epsilon,b[) \Big)^{n-1}\\
&=& \int_{[b,\infty[} e^{-\theta x}\nu_a(dx)
\frac{1}{\nu_a[\epsilon,\infty[} \sum_{n=1}^\infty
\Big(\frac{\int_{[\epsilon,b[} e^{-\theta
    x}\nu_a(dx)}{\nu_a[\epsilon,\infty[}  \Big)^{n-1}\\ 
&=& \frac{\int_{[b,\infty[} e^{-\theta
    x}\nu_a(dx)}{\nu_a[\epsilon,\infty[-\int_{[\epsilon,b[} e^{-\theta
    x}\nu_a(dx)}. 
\end{eqnarray*}

Par stationnarité des accroissements de $B$, la loi de $\mc_{a,b}$ est
invariante par translation si bien que la suite $(T'_n-a)_{n \ge 1}$
a même loi que la suite $(T_n)_{n \ge 1}$ où $T_1<T_2<\ldots$ sont les
instants successifs de $\mc_{a,b} \cap \rrf_+^*$. On peut donc énoncer
le résultat suivant. 

\begin{prop}~\label{lien entre mesure de Levy et durees successives}
Soient $a \ge b > 0$. Notons $T_1<T_2<\ldots$ les instants successifs
de $\mc_{a,b} \cap \rrf_+^*$. Les durées $T_2-T_1,T_3-T_2,\ldots$ sont
indépendantes entre elles et avec $T_1$, de même loi ; leur transformée de
  Laplace est donnée par $$\eef[e^{-\theta (T_2-T_1)}] =
  \frac{\int_{[b,\infty[} e^{-\theta x} \nu_a(dx)}{\int_{\rrf_+^*} (1
    - \iif_{[x  < b]}e^{-\theta x}) \nu_a(dx)}$$
où $\nu_a$ est la mesure de Lévy du subordinateur $\sigma^a$.
\end{prop}

{\parindent 0cm  {\bf Remarque : une autre preuve du théorème~\ref{Loi
      des durees}}}.  

Nous venons de redémontrer l'indépendance des variables
$T_1,T_2-T_1,T_3-T_2,\ldots$ et l'équidistribution des durées
$T_n-T_{n-1}$. Une fois ce point établi, on retrouve facilement la
loi $\mu_{a,b}$ de $T_1$ et la loi $\nu_{a,b}$ de $T_2-T_1$ à partir
du théorème~\ref{fonction de correlation} : pour tout $t>s>0$,
$$\frac{ds}{\pi\sqrt{ab}} = \eef[N_{a,b}(ds)] = \sum_{m=1}^\infty
P[T_m \in ds] = \sum_{m=1}^\infty (\mu_{a,b}*\nu_{a,b}^{*(m-1)})(ds)$$
et
\begin{eqnarray*} 
\frac{ds}{\pi\sqrt{ab}} h_{a,b}(t-s) dt = \eef[N_{a,b}(ds)N_{a,b}(dt)] 
&=& \sum_{m=1}^\infty \sum_{n=1}^\infty P[T_m \in ds\ ;\ T_{m+n} \in dt]\\ 
&=& \sum_{m=1}^\infty \sum_{n=1}^\infty
(\mu_{a,b}*\nu_{a,b}^{*(m-1)})(ds)\ \nu_{a,b}^{*n}(dt-s)
\end{eqnarray*}
d'où  
$$h_{a,b}(r) dr = \sum_{n=1}^\infty \nu_{a,b}^{*n}(dr)$$
On en déduit les égalités suivantes pour les transformées de Laplace 
$$\frac{1}{\pi\sqrt{ab}\theta} = \frac{\lc \mu_{a,b}(\theta)}{1-\lc
  \nu_{a,b}(\theta)}$$ 
$$\lc h_{a,b}(\theta) = \frac{\lc \nu_{a,b}(\theta)}{1-\lc \nu_{a,b}(\theta)}$$
On retrouve ainsi les lois $\nu_{a,b}$ et $\mu_{a,b}$ décrites dans le
théorème~\ref{Loi des durees} par l'intermédiaire de leur
transformée de Laplace.

\section{Description de la mesure de Lévy du subordinateur $\sigma^a$}

Dans toute cette partie, on étudie la fonction $G_a$ définie par
$G_a(r) = \nu_a]r,\infty[$ pour $r>0$ et $G_a(r)=0$ pour $r \le 0$. 

\subsection{Equations de convolution vérifiées par $G_a$}

En utilisant les égalités ci-dessus et la proposition~\ref{lien entre
mesure de Levy et durees successives}, on obtient lorsque $a \ge b >
0$, 
\begin{equation*}
\lc h_{a,b}(\theta) = \frac{\lc \nu_{a,b}(\theta)}{1-\lc
  \nu_{a,b}(\theta)} = \frac{\int_{[b,\infty[} e^{-\theta x} 
  \nu_a(dx)}{\int_{\rrf_+^*} (1 - e^{-\theta x}) \nu_a(dx)}.
\end{equation*}
Comme 
$$\int_{\rrf_+^*} (1 - e^{-\theta x})\ \nu_a(dx) = \int_0^\infty \theta\
e^{-\theta r}\ G_a(r)\ dr = \theta\ \lc G_a (\theta),$$
l'égalité précédente s'écrit aussi 
$$\lc h_{a,b}(\theta)\ \lc G_a (\theta) = \frac{1}{\theta}
\int_{[b,\infty[} e^{-\theta x}\ \nu_a(dx),$$ 
soit en termes de convolution : 
$$h_{a,b} * G_a = \iif_{\rrf_+} * (\iif_{[b,\infty[} \nu_a).$$
Pour tout $x \ge 0$, on a donc 
\begin{equation}~\label{lien entre $h_{a,b}$ et mesure de Levy}
  \int_0^x h_{a,b}(y)\ G_a(x-y)\ dy = G_a(b) - G_a(x \vee b).
\end{equation} 
Comme $h_{a,b}(y)=0$ pour $y \le b$, cette égalité donne la valeur de $G_a(x)$
pour $x \ge b$ en fonction de $G_a(b)$ et des valeurs $G_a(z)$ pour $z
\in [0,x-b]$. Comme $G_a(x)=\sqrt{2/\pi x}$ pour tout $x \in [0,a]$,
$G_a$ est connue {\it a fortiori} sur $[0,b]$, et la relation
~(\ref{lien entre $h_{a,b}$ et mesure de Levy}) permet de calculer
par récurrence la valeur de $G_a$ sur les intervalles de la forme
$[nb,(n+1)b]$ avec $n \in \nnf$. Les formules deviennent vite
compliquées même si elles se simplifient un peu dans le cas où $b=a$. 

%En effet, comme  pour tout $x \in [0,a]$, 
%$$G_a(x) = \sqrt{\frac{2}{\pi x}}$$
%on obtient pour tout $x \in [a,2a]$, 
%$$G_a(x) = G_a(a) - \int_a^x \frac{2}{\pi y} \sqrt{\frac{y-a}{a}} G_a
%  (x-y) dy = \sqrt{\frac{2}{\pi a}}\Big(1 - \frac{2}{\pi} \int_a^x 
%  \sqrt{\frac{y-a}{x-y}}\ \frac{dy}{y} \Big)$$

Un autre cas particulier plus intéressant est le cas limite où $b \to
0$. En effet, pour tout $b \le a$, $\sqrt{b}G_a(b) = \sqrt{2 \pi}$ et
pour tout $y>0$, $\pi \sqrt{b} h_{a,b}(r) \to (y \wedge a)^{-1/2}$
quand $b \to 0$. En multipliant par $\pi \sqrt{b}$ l'égalité~\ref{lien
  entre $h_{a,b}$ et mesure de Levy} et en faisant tendre $b$ vers
$0$, on obtient donc par convergence dominée 
\begin{equation}~\label{equation de convolution}
\int_0^x (y \wedge a)^{-1/2}\ G_a(x-y)\ dy = \sqrt{2 \pi} 
\end{equation}
Cette égalité fournit une expression simple de la transformée de
Laplace de $G_a$.  

\subsection{Expressions de la transformée de Laplace de $G_a$ et applications}
 
\begin{theo}~\label{Transformee de Laplace de $G_a$}
La transformée de Laplace de $G_a$ est donnée par 
$$\lc G_a(\theta) = \sqrt{2 \pi a}
M \big(-\frac{1}{2};\frac{1}{2};-\theta a \big)^{-1}$$
où $M(a,b,\cdot) = {}_1F_1(a,b,\cdot)$ est la fonction hypergéométrique
définie par : 
$$M(a,b,x) = \sum_{n=0}^{+\infty} \frac{(a)_n}{(b)_n}\frac{x^n}{n!},$$
en notant $(a)_n = a (a+1) \cdots (a+n-1)$. 
%On vérifie immédiatement que
%$$\frac{\partial}{\partial x}\ {}_1F_1(a,b,x) = \sum_{n=1}^{+\infty}
%\frac{a_{n+1}}{b_{n+1}}\frac{x^n}{n!} = \frac{a}{b}\ {}_1F_1(a+1,b+1,x).$$
\end{theo}

\begin{demo}
L'équation~\ref{equation de convolution} ci-dessus entraîne
imédiatement l'égalité
$$\lc G_a(\theta) \times \int_0^\infty (r \wedge a)^{-1/2}
  \ e^{-\theta r} dr = \frac{\sqrt{2 \pi}}{\theta}.$$
Mais 
\begin{eqnarray*}
\int_0^\infty (r \wedge a)^{-1/2} \ e^{-\theta r} dr 
&=& \int_a^\infty a^{-1/2}\ e^{-\theta r} dr\ +\ \int_0^a r^{-1/2}\
 e^{-\theta r} dr\\ 
&=& \frac{a^{-1/2}}{\theta} e^{-\theta a}\ +\ \sum_{n=0}^{+\infty}
\frac{(-\theta)^n}{n!}\ \frac{a^{n+1/2}}{n+1/2}\\ 
&=& \frac{a^{-1/2}}{\theta} \Big(\sum_{n=0}^{+\infty} \frac{(-\theta
  a)^n}{n!}\ -\ \sum_{n=1}^{+\infty} \frac{(-\theta
  a)^n}{(n-1)!(n-1/2)}\Big)\\ 
&=& \frac{a^{-1/2}}{\theta} \sum_{n=0}^{+\infty} \frac{(-\theta
  a)^n}{n!} \Big(1\ -\ \frac{n}{n-1/2}\Big)\\
&=& \frac{a^{-1/2}}{\theta} M \big(-\frac{1}{2};\frac{1}{2};-\theta a \big),
\end{eqnarray*}
puisque pour tout $n \in \nnf$, 
$$1\ -\ \frac{n}{n-1/2} = \frac{-1/2}{n-1/2} =
\frac{(-1/2)_n}{(1/2)_n}.$$
ce qui démontre le théorème \hfill $\square$
\end{demo}

La transformée de Laplace fournit des indications sur le
comportement asymptotique de la queue de $G_a$. En effet, 
par prolongement analytique, on déduit du théorème~\ref{Transformee de
  Laplace de $G_a$} l'égalité
$$\int_0^\infty e^{\theta r}\ G_a(r)\ dr = \sqrt{2 \pi a}\ M \big(-\frac{1}{2};\frac{1}{2};\theta a \big)^{-1}$$
pour tout réel $\theta$ tel que $\theta a < \rho$, où $\rho$ 
est l'unique zéro dans $\rrf_+$ de la fonction
$M(-\frac{1}{2};\frac{1}{2};\cdot)$. De plus, quand $\varepsilon \to 0+$, 
$$\int_0^\infty e^{-\varepsilon r} e^{\rho r/a}\ G_a(r)\ dr = \sqrt{2 \pi
  a}\ M \big(-\frac{1}{2};\frac{1}{2};\rho - a \varepsilon
\big)^{-1} \sim \frac{\sqrt{2 \pi a}}{\lambda a \varepsilon}.$$
avec $\lambda = M(\frac{1}{2};\frac{3}{2};\rho)$. A l'aide du théorème 
taubérien de Hardy ou de Karamata (voir \cite{Feller} au chapitre
XIII), on en déduit le résultat suivant. 

\begin{prop} Notons $\rho$ l'unique zéro dans $\rrf_+$ de la fonction
$M(-\frac{1}{2};\frac{1}{2};\cdot)$ et
$\lambda=\frac{\partial}{\partial x}\
M(-\frac{1}{2};\frac{1}{2};\cdot) (\rho) = M(\frac{1}{2};\frac{3}{2};\rho)$. Alors quand $x \to
+\infty$, 
$$\int_0^x e^{\rho r/a}\ G_a(r)\ dr \sim \frac{\sqrt{2
    \pi}}{\lambda\sqrt{a}}\ x.$$
\end{prop}

On peut s'attendre à ce que la queue $G_a$ soit suffisamment régulière
pour que $e^{\rho r/a}\ G_a(r)$ ait une limite quand $r \to
+\infty$, ce qui conduit à la conjecture suivante
sur le comportement asymptotique de la queue de $\nu_a$. 

\paragraph{Conjecture} Avec les mêmes notations que ci-dessus, 
$$G_a(r) \sim \frac{\sqrt{2 \pi}}{\lambda\sqrt{a}}\ e^{-\rho r/a}.$$

Cette conjecture est confortée par la représentation graphique de $\ln
G_1$ donnée à la fin de l'article.

De même, si $a \ge b > 0$ et si $T_1<T_2<\ldots$ sont
les instants successifs de $\mc_{a,b} \cap \rrf_+^*$, on peut
conjecturer de même que la variable aléatoire $T_2 - T_1$ a une queue
(et une densité) à décroissance exponentielle, puisque d'après la
proposition~\ref{lien entre mesure de Levy et durees successives}, on
a pour tout $\theta>0$, 
$$1 - \eef[e^{-\theta (T_2-T_1)}] =  
  \frac{\int_{[0,\infty[} (1-e^{-\theta x}) \nu_a(dx)}{\int_{\rrf_+^*} (1
    - \iif_{[x  < b]}e^{-\theta x}) \nu_a(dx)} = 
  \frac{\sqrt{2\pi}\ \theta \lc G_a(\theta)}{2b^{-1/2} + \int_{[0,b]}
    (1-e^{-\theta x}) x^{-3/2}dx}$$
Cette formule s'étend par prolongement analytique à tout $\theta >
-\rho/a$, et il existe une constante $C \in \rrf_+^*$ telle que
$$\eef[e^{(\frac{\rho}{a}-\varepsilon)(T_2-T_1)}] \sim C/\varepsilon\ \text{
  quand }\ \epsilon \to 0+.$$

\subsection{Expression de $G_a$ sous forme de série}

La transformée de Laplace de $G_a$ peut être réécrite comme suit. 

\begin{eqnarray*}
\lc G_a(\theta) &=& \frac{\sqrt{2\pi}}{\theta}\ \Big(\int_0^\infty
  r^{-1/2}\ e^{-\theta r} dr + \int_a^\infty (a^{-1/2} - r^{-1/2})\
  e^{-\theta r} dr \Big)^{-1}\\ 
&=& \frac{\sqrt{2\pi}}{\theta}\ \Big(\sqrt{\frac{\pi}{\theta}} +
\frac{1}{2\theta} \int_a^\infty r^{-3/2}e^{-\theta r} dr \Big)^{-1}\\ 
&=& \sqrt{\frac{2}{\theta}}\ \Big(1 + \frac{1}{2 \sqrt{\pi\theta}}
  \int_a^\infty r^{-3/2}e^{-\theta r} dr \Big)^{-1}. 
\end{eqnarray*}
Pour tout $\theta$ suffisamment grand, on a donc 
$$\lc G_a(\theta) = \sqrt{\frac{2}{\theta}}\ \sum_{n=0}^{+\infty} (-1)^n
\Big(\frac{1}{2 \sqrt{\pi\theta}} \int_a^\infty r^{-3/2}e^{-\theta r}
dr \Big)^n.$$
Remarquons que l'application 
$$\theta \mapsto \Big(\frac{1}{2 \sqrt{\pi\theta}} \int_a^\infty
r^{-3/2}e^{-\theta r} dr \Big)^n$$ 
est la transformée de Laplace du produit de convolution de $r \mapsto
(2\pi)^{-1}\iif_{[r>0]}r^{-1/2}$ avec $r \mapsto \iif_{[r>a]}r^{-3/2}$.
Mais pour $r>a$,  
\begin{eqnarray*}
\int_a^r \frac{dy}{(r-y)^{1/2}y^{3/2}} = \int_a^r
\sqrt{\frac{y}{r-y}}\ \frac{dy}{y^2} = \int_{a/(r-a)}^\infty
\sqrt{z}\ \frac{dz}{rz^2} = \frac{2}{r}\sqrt{\frac{r-a}{a}},
\end{eqnarray*}
grâce au changement de variable $z = y/(r-y)$ d'où $1/z = r/y - 1$. 
On déduit des calculs précédents une expression de $G_a$.   

\begin{theo}~\label{Calcul de $G_a$}
Notons $u$ et $h_{\infty,a}$ 
les applications définies par 
$$u(r) = \iif_{[r>0]} \sqrt{\frac{2}{\pi r}}\ \text{ et }\
h_{\infty,a}(r) = \frac{\iif_{[r>a]}}{\pi r}\sqrt{\frac{r-a}{a}}.$$
Alors
$$G_a = \sum_{n=0}^{+\infty} (-1)^n (u * h_{\infty,a}^{*n})$$
\end{theo}

{\bf Remarque. }Dans cet énoncé, la notation $h_{\infty,a}$ a été
choisie par analogie avec les 
notations du théorème~\ref{fonction de correlation}. Cependant nous
n'avons pas d'explication à l'égalité formelle $G_a = u -
u*g_{\infty,a}$, dans laquelle la convolution par $u$ s'apparente à
une intégration fractionnaire. 

{\bf Calcul de $G_a$.} Pour tout $n \in \nnf$, $u * h_{\infty,a}^{*n}$
est à support dans $[na,+\infty[$. Pour calculer $G_a$ sur
un intervalle $[0,Na]$ avec $N \in \nnf$, il suffit donc de faire
varier $n$ de $0$ à $N-1$ dans la série qui intervient dans le
théorème~\ref{Calcul de $G_a$}. Par exemple, pour tout $r \in [0,2a]$,
$G_a(r) = u(r) - u * h_{\infty,a}(r)$. 

Or pour tout $r \in [a,+\infty[$,  
$$u * h_{\infty,a}(r) = \sqrt{\frac{2}{\pi a}}\ \frac{1}{\pi} \int_a^r
\sqrt{\frac{x-a}{r-x}}\ \frac{dx}{x}.$$
Le changement de variable 
$$t = \sqrt{\frac{x-a}{r-x}},\ \text{ soit }\ x =
\frac{rt^2+a}{t^2+1},$$
$$\frac{dx}{x} = \Big(\frac{2rt}{rt^2+a}- \frac{2t}{t^2+1}\Big)dt =
\Big(\frac{1}{t^2+1} - \frac{a}{rt^2+a} \Big) \frac{2dt}{t},$$
montre que pour tout $r \in [a,+\infty[$,
$$u * h_{\infty,a}(r) = \sqrt{\frac{2}{\pi a}}\ \frac{2}{\pi}
\int_0^\infty \Big(\frac{1}{t^2+1} - \frac{a}{rt^2+a} \Big)\ dt 
= \sqrt{\frac{2}{\pi a}}\ \Big(1 - \sqrt{\frac{a}{r}} \Big)
= \sqrt{\frac{2}{\pi a}} - \sqrt{\frac{2}{\pi r}}.$$ 

Pour tout $r \in [a,2a]$, on a ainsi
$$G_a(r) = u(r) - (u * h_{\infty,a})(r) = 2 \sqrt{\frac{2}{\pi r}} -
\sqrt{\frac{2}{\pi a}}.$$

\subsection{Représentations graphiques de $G_1$ et $\ln G_1$}

La figure ci-dessous montre le graphe de $G_1$, calculé à l'aide de
Scilab. Je remercie J.M. Decauwert pour son aide précieuse en la matière. 
Notons que $G_a$ se déduit de $G_1$ par changement d'échelle. Plus
précisément, la relation
$G_a(r) = G_1(r/a)/\sqrt{a}$ découle de l'égalité $\lc G_a(\theta) =
\sqrt{a} \lc G_1(a\theta)$ par injectivité de la transformation de Laplace.

\begin{figure}[hbtp]
  \begin{center} 
\includegraphics[width=16cm]{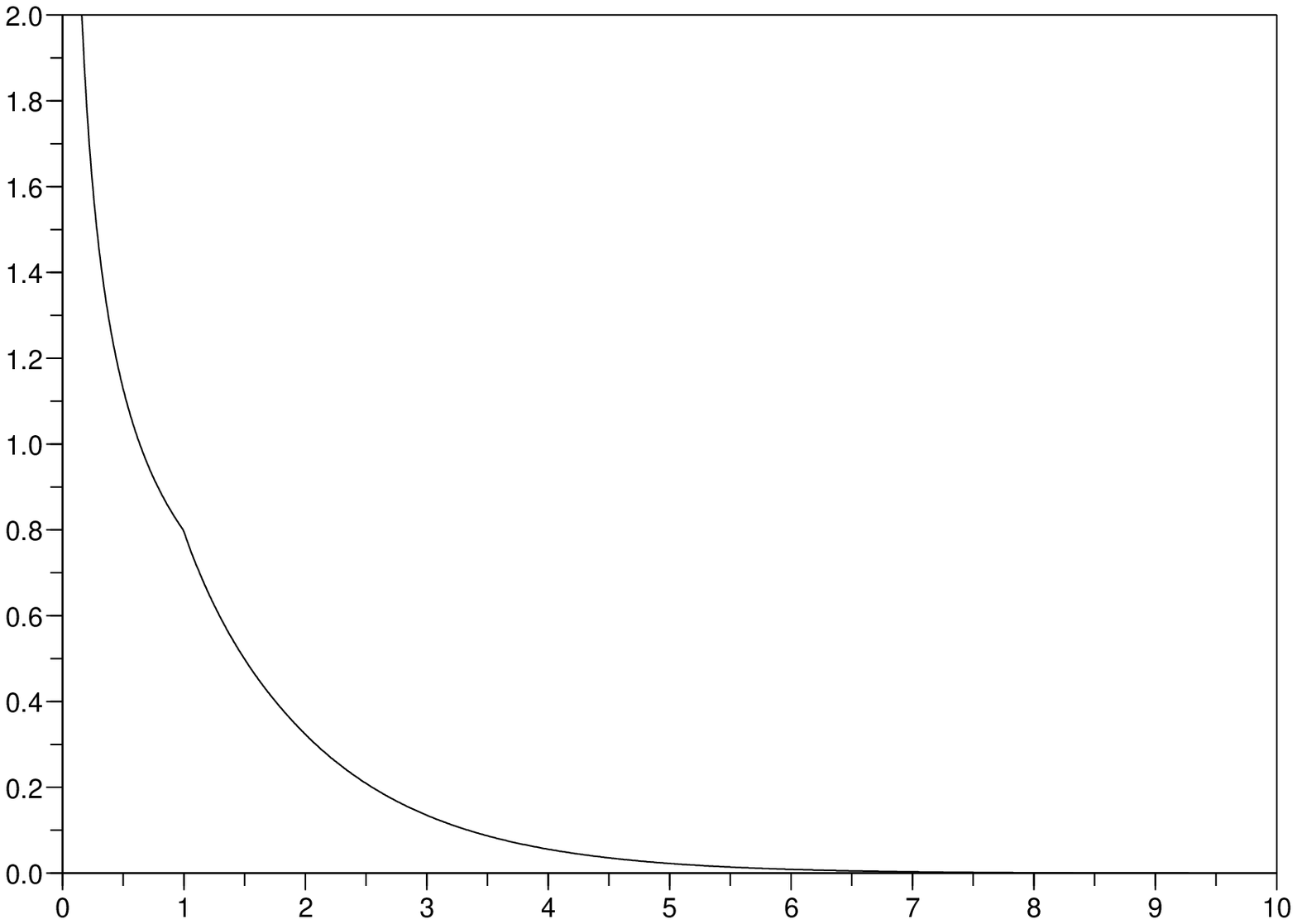}
  \end{center}                                    
  Figure 1. --- Représentation graphique de $G_1$.               
\end{figure}

La conjecture sur la croissance exponentielle de $G_a$
faite plus haut peut se réécrire sous la forme $\ln G_1(r) + \rho
r \to \ln (\sqrt{2\pi}/\lambda)$ quand $r \to +\infty$. Il est donc
intéressant de dessiner le 
graphe de $\ln G_1$. La figure ci-dessous suggère en effet que la
courbe représentant $\ln G_1$ possède une asymptote de pente
voisine de -0,9. Le graphe a été limité à l'intervalle $[0,5]$ en
raison de problèmes d'instabilité numérique qui apparaissent au-delà. 

\begin{figure}[hbtp]
  \begin{center} 
\includegraphics[width=16cm]{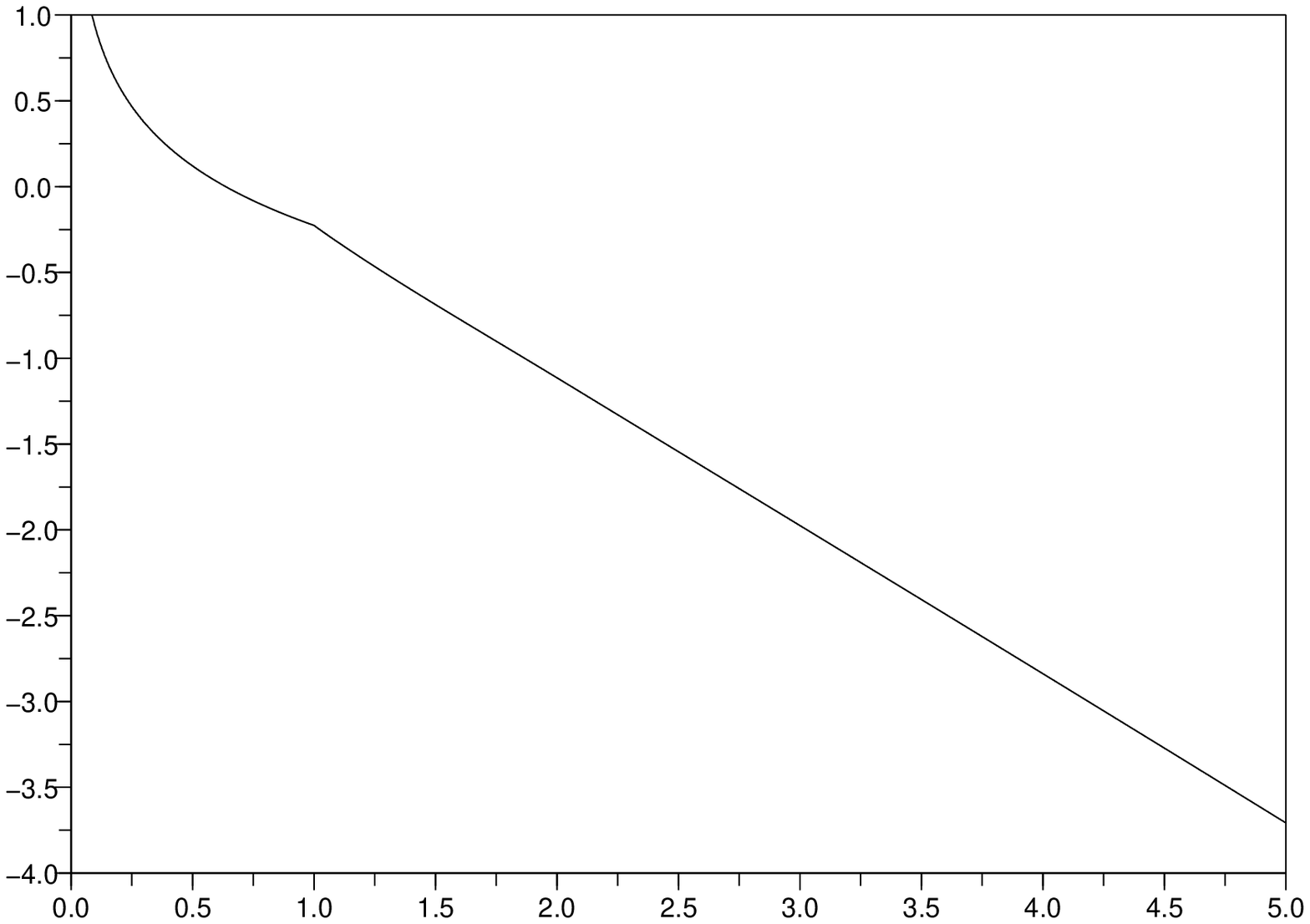}    
%    \leavevmode \input{graphe_G_1.eps} %\label{fig:}
  \end{center}                                    
  Figure 1. --- Représentation graphique de $\ln G_1$            
\end{figure}

\vfill
\begin{flushleft}
%adresse de l'auteur
Christophe LEURIDAN\\
INSTITUT FOURIER\\
Laboratoire de Math\'ematiques\\
UMR5582 (UJF-CNRS)\\
BP 74\\
38402 St MARTIN D'H\`ERES Cedex (France)\\
Christophe.Leuridan@ujf-grenoble.fr
\end{flushleft}

\end{document}